# Une méthode pour obtenir la fonction génératrice algébrique d'une série.

**From FPSAC, Formal Power Series and Algebraic Combinatorics, Florence, june 1993.**

Simon Plouffe
LACIM
Université du Québec à Montréal
Mars 1993

**Résumé**

Nous décrivons ici une méthode expérimentale permettant de calculer de bons candidats pour une forme close de fonctions génératrices à partir des premiers termes d'une suite de nombres rationnels. La méthode est basée sur l'algorithme LLL[1] et utilise deux programmes de calcul symbolique, soit MapleV et Pari-GP. Quelques résultats sont présentés en appendice. Cette méthode a été testée sur toute la table de suites du livre , *The New book of Integer Sequences*, de N.J.A Sloane et S. Plouffe (en préparation). Ainsi, nous avons obtenu de cette façon la fonction génératrice,

$$\frac{z + (z+1)^{1/2}(1-3z)^{1/2} - 1}{2(z^2(z+1)^{1/2}(1-3z)^{1/2})}$$

pour la suite: 1, 2, 6, 16, 45, 126, 357, 1016, 2907, 8350, 24068, 69576, 201643, 585690,... qui apparaît en page 78 du livre de Louis Comtet, *Adanced Combinatorics*.

---

[1] Nommé ainsi à cause des travaux de Lenstra, Lenstra et Lovasz.

# INTRODUCTION

## Les suites P-récurrentes

On dit qu'une suite $a_n$ d'entiers ou de nombres rationnels est P-récurrente si on a une récurrence de la forme :

(1) $$a_n P_0(n) = a_{n-1} P_1(n) + a_{n-2} P_2(n) + \ldots + a_{n-k} P_k(n),$$

où les $P_i(n)$, $0 \leq i \leq k$, sont des polynômes à coefficients rationnels. D'un autre point de vue, (1) équivaut à dire que la fonction génératrice

$$A(x) = \sum_{n \geq 0} a(n) x^n$$

satisfait à une équation différentielle linéaire à coefficients polynômiaux. Il n'est donc pas possible d'espérer pouvoir trouver systématiquement une forme close pour $A(x)$. Si l'on connait N termes d'une suite $a_n$ satisfaisant (1) avec $N > \sum (\deg(P_i)+1)$, on peut déterminer les $P_i(x)$ par une approche de coefficients indéterminés.

Le programme [gfun] de la librairie "share" de MapleV procède ainsi via la commande "*listtorec*" pour obtenir une équation de la forme (1) pour une suite dont les N premiers termes sont :$[a_0, a_1, a_2, \ldots, a_N]$. Par exemple, la suite #1173 de [Sl,Tutte], relative aux cartes planaires, est une suite P-récurrente et on obtient:

## Exemple 1

```
> suite:=[1, 1, 0, 1, 3, 12, 52, 241, 1173, 5929, 30880, 164796, 897380,
4970296, 27930828, 158935761, 914325657, 5310702819, 31110146416,
183634501753, 1091371140915]:

> recsuite:=listtorec(suite):

                    3              /            59   2          3\
{(1/4 + 7/8 n - 9/8 n ) a(n) + |- 5/4 + 2/3 n + ---- n  - 13/3 n | a(n+1)
                                   \            12               /
                 2       3
      + (- 1 - 2/3 n + n  + 2/3 n ) a(n+2), a(1) = 1, a(2) = 1}.
```

On peut automatiquement traduire cette P-récurrence en programme avec la commande "rectoproc",

```
> Sl1173:=rectoproc(recsuite,a(n));

Sl1173 := proc(n) options remember;
if not type(n,nonnegint) then ERROR(`invalid arguments`) fi;
1/8*(27*procname(-2+n)*n^2+104*procname(n-1)*n^2-81*
procname(-2+n)*n-430*procname(n-1)*n+60*
procname(2+n)+414*procname(n-1))/(2*n^2-3*n+1)
end;
```

Notons que la procédure fournie permet de calculer autant de termes que l'on veut en temps linéaire par rapport à n. Cela nous servira par la suite. On dit de la P-récurrence précédente qu'elle est est de type (2,3), c'est à dire de degré 2 et à 3 termes. Lorsque la P-récurrence obtenue est de type (d,2) ou d est le degré, le rapport de deux termes successifs est une fraction rationnelle. On peut obtenir une expression hypergéométrique pour le terme général de la suite $a_n$.

Dans un même ordre d'idée, on peut s'intéresser au cas où la fonction génératrice A(z) est algébrique, i.e.:

$$\sum_{j=0}^{n} P_j(z)A(z)^j = 0$$

pour certains polynômes $P_j(z)$. Ce que nous proposons ici est une méthode expérimentale pour passer de la P-récurrence trouvée via [gfun] à l'équation algébrique.

Bien entendu on pourrait toujours tenter de résoudre le système

$$\sum_{0 \le j,k \le m} c_{j,k} S(z)^j z^k = 0$$

satisfait par la fonction génératrice, où S(z) est la série génératrice des $a_n$ et les $c_{j,k}$ sont les indéterminées. C'est cette approche qu'utilise la procédure "*listtoalgeq*" de [gfun]. Ainsi pour la suite des nombres de Catalan, voir [Sl], [Comtet]

**Exemple 2**

```
> suite := [1, 1, 2, 5, 14, 42, 132, 429, 1430, 4862, 16796, 58786, 208012,
742900, 2674440, 9694845, 35357670, 129644790, 477638700, 1767263190,
6564120420];

> listtoalgeq(suite,S(z));

                                   2
                  1 - S(z) + z S(z)
```

En résolvant par rapport à S(z) et en prenant la racine positive on obtient,

```
> solve(",S(z));

                                        1/2
                       1 - (1 - 4 z)
              1/2  ---------------
                         z
```

Cependant cette approche est limitée pour deux raisons principales. D'abord le nombre d'équations et d'inconnues explose rapidement; de plus, comme les calculs s'effectuent en précision infinie et à l'aide d'un interprète, la résolution d'un gros système devient vite coûteuse en temps et en espace. L'expérience suggère qu'il est plus facile d'obtenir une P-récurrence explicite qu'une équation algébrique explicite. On peut donc penser rechercher d'abord une P-récurrence pour ensuite obtenir par un processus efficace une équation algébrique. C'est précisément ce genre d'approche que nous allons décrire.

## 2. L'algorithme LLL.

L'algorithme LLL permet, étant donné une base d'un espace vectoriel de dimension finie qui génère un réseau, de trouver une base réduite de ce même espace. Si le réseau, engendré par les vecteurs de la base, est formé de vecteurs à coordonnées entières, celui-ci trouvera en temps polynômial, des vecteurs plus courts qui généreront le même réseau. L'algorithme existe également dans la version qui permet d'avoir des vecteurs à coordonnées rationnelles. On peut se référer aux articles originaux [LLL] et [Kannan]. Il permet, entre autres, de trouver numériquement un polynôme dont un nombre "réel" en virgule flottante est la racine.

Cet algorithme est implanté depuis la version V de Maple et porte le nom de "minpoly". Il fait appel à l'algorithme LLL. Disons simplement qu'il permet de résoudre, entre autres, numériquement le problème exactement inverse de trouver une racine d'un polynôme.

Le problème inverse étant : si un nombre réel est donné, que l'on suppose algébrique, de quel polynôme minimal est-il racine ? Mentionnons dès maintenant qu'on parle ici d'un nombre réel donné avec une certaine précision numérique. Ce nombre peut se convertir en nombre rationnel équivalent à la précision numérique voulue. On ne pourra (une fois l'opération réussie) qu'isoler un polynôme qui *semble* avoir ce nombre réel comme racine. Il serait un peu long de donner tous les détails qui font qu'aujourd'hui ce problème est pour ainsi dire *numériquement résolu.*

Mentionnons qu'au moins trois programmes de calcul symbolique ont implanté cet algorithme: soit Maple, Mathematica et Pari-GP. La meilleure implémentation de cet algorithme et la plus rapide nous semble être est celle de Pari-GP [Pari].

Cette procédure accepte donc en entrée un nombre décimal tronqué et donne (selon la précision numérique en vigueur) un polynôme **minimal** P(x) dont serait racine. La précision numérique en vigueur est celle que l'utilisateur demande. En pratique, elle est d'environ 500 chiffres décimaux avec Pari-GP. Le degré maximal du polynôme que l'on puisse demander dépend largement de cette précision. En pratique, la limite est un polynôme de degré 20. Ceci est tout de même suffisant pour obtenir des résultats intéressants.

Si la fonction génératrice S(z) correspondant à la suite est algébrique et si z=1/m, avec m entier, S(1/m) sera un nombre algébrique. C'est précisément ici que l'on utilise l'algorithme LLL. Avec la procédure "*rectoproc*" , on peut obtenir des milliers de termes qui serviront à évaluer S(z) en un point 1/m "très petit". Ainsi le résultat sera un nombre algébrique approché à une grande précision numérique. On calcule alors le polynôme dont S(1/m) est racine. On peut vérifier la vraisemblance du résultat en répétant ce calcul pour S(1/(m+1)), S(1/(m+2)), S(1/(m+3)), ... . Il se trouve que la version de LLL du programme Pari-GP est extrêmement efficace. Non seulement la procédure (qui s'appelle "*algdep*") retourne en général le bon polynôme, mais de surcroît il est simplifié au maximum. De plus, les solutions trouvées sont *numériquement stables*; elles sont stables au point qu'elles permettent de reconstruire la fonction génératrice algébrique.

Voici donc l'algorithme correspondant. On utilise ici un interface qui permet de passer de Maple à Pari-Gp et vice versa dans la même session. Les commandes "listtorec,rectoproc listtoseries" font partie du programme [gfun], la commande ALGDEP du programme Pari-Gp; le reste des commandes, comme (evalf,interp,solve) font partie de MapleV en version standard. La commande "*interp*" de Maple permet simplement de calculer le polynôme d'interpolation de Newton et utilise les différences finies.

**1)** listtorec(suite);
**2)** rec:=rectoproc(suite);
**3)** ser:=listtoseries(suite,z,ogf);
**4)** for i=1 from 1 to **nombre** do
   v(i):=subs(z=1/(**m**+i-1),ser);
   vf(i):=evalf(v(i),precision);
   polynôme(i):=**PARI**(**ALGDEP**(vf(i),**degré**));
  od;
**5)** eqalgébrique:=interp(polynôme(i),t,**m**);
**6)** alendroit:=subs(t=1/z,eqalgébrique);
**7)** solve(alendroit,t);

Illustrons cet algorithme en donnant un exemple. Nous prendrons une suite qui apparaît dans le Journal of Combinatorial Theory B, vol. 21 (1976) pp. 71-75, et qui est relative aux tournois. Article de J.W. Moon. Ici nous prendrons **nombre**=13, **degré**=2 et **m**=100. C'est à dire que l'on suppose ici que le degré de l'équation est 2, mais les coefficients sont de degré jusqu'à 12.

**Exemple 3**

```
> suite:=[1, 1, 1, 3, 16, 75, 309, 1183, 4360, 15783, 56750, 203929, 734722,
2658071,  9662093,  35292151,  129513736,  477376575,  1766738922, 6563071865,
24464169890];

> listtorec(suite,a(n));

[{a(2) = 1, a(1) =
                        2                           2
     (- 2/3 n - 4/3 n ) a(n) + (- 1 + n + n ) a(n + 1)

                               2                          2         3
         + (1/2 - 1/3 n - 1/6 n ) a(n + 2) + 1/2 + 1/6 n - 1/6 n  + 1/2 n }]

> rec:=rectoproc(",a(n));

rec :=proc(n)
 options remember;
     if not type(n,nonnegint) then ERROR(`invalid arguments`) fi;
     (28*procname(n-2)*n-24*procname(n-2)-8*
     procname(n-2)*n^2+6*procname(n-1)-18*procname(n-1)*n+6*
     procname(n-1)*n^2-27+41*n-19*n^2+3*n^3)/(-3-2*n+n^2)
 end;
```

Voici la première valeur décimale, c'est la série S(z) évaluée en z=1/100. Le lecteur attentif reconnaîtra les premiers termes de la suite :1,1,3,16,... dans le développement décimal du nombre.

```
vf(1)=1.010103167821282371655205556160928600562159888369694333057529335554251
```

5029460058952354762187795026581944514416380788705715044395043768728954722738516149864952340103813169557832245178542753139285380720304392389878530808969233130
46663

Voici les polynômes trouvés par **ALGDEP** de Pari-Gp en quelques secondes de calcul seulement.

$$922556408004\,x^2 - 9041033588479200\,x + 9131435376040000$$
$$980100000000\,x^2 - 9799999702020000\,x + 9897020403050401$$
$$1040604010000\,x^2 - 10614139675759200\,x + 10718190400203216$$
$$1104189046416\,x^2 - 11486856353906376\,x + 11598369273824917$$
$$1170979365924\,x^2 - 12421725705345216\,x + 12541154909460736$$
$$1241102946304\,x^2 - 13422503799519360\,x + 13550326173504225$$
$$1314691560000\,x^2 - 14493133991044800\,x + 14629850124065296$$
$$1391880848400\,x^2 - 15637754317171560\,x + 15783889435204501$$
$$1472810396836\,x^2 - 16860705112257696\,x + 17016810038701632$$
$$1557623810304\,x^2 - 18166536843458976\,x + 18333188987567041$$
$$1646468789904\,x^2 - 19560018171877920\,x + 19737822545544400$$
$$1739497210000\,x^2 - 21046144243456200\,x + 21235734506893941$$

```
> interp(polynôme(i),t,100);
```

$$(1 - 9z + 32z^2 - 57z^3 + 54z^4 - 24z^5 + 4z^6 - t + 10tz - 42tz^2$$
$$+ 98tz^3 - 137tz^4 + 112tz^5 - 48tz^6 + 8tz^7 + t^2z - 8t^2z^2$$
$$+ 26t^2z^4 - 44t^2z^5 + 41t^2z^6 - 20t^2z^7 + 4t^2z^8) / z^8$$

C'est l'équation algébrique recherchée; elle est de degré 2 en la variable t. Il ne reste qu'à résoudre cette équation et l'une des deux solutions (la positive) est,

$$\frac{-1/2\,(-1 + 10z - 42z^2 + 98z^3 - 137z^4 + 112z^5 - 48z^6 + 8z^7)}{(z^2\,(2z-1)^2\,(z-1)^4)}$$

$$\frac{(-(-1+4z)\,(2z-1)^4\,(z-1)^8)^{1/2})}{(z^2\,(2z-1)^2\,(z-1)^4)}$$

En développant en série de Taylor cette fonction génératrice, on retrouve bien notre suite de départ. On en conclut que c'est *probablement* la fonction génératrice algébrique de cette suite.

**Conclusion** : notre méthode permet donc de trouver souvent une fonction génératrice algébrique de degré élevé en autant que l'on puisse trouver une P-récurrence. Les temps de calcul sont relativement peu élevés puisque ceux-ci ont été vérifiés avec un micro-ordinateur Macintosh SE/30. (4 mips). Dans [Plo], plus de 32 fonctions génératrices algébriques ont été trouvées grâce à cette méthode. Nous en présentons quelques unes en appendice.

1, 2, 9, 54, 378, 2916, 24057, 208494, 1876446, 17399772, 165297834, 1602117468, 15792300756, 157923007560, 1598970451545, 16365932856990
**Réf. : CJM 15 254 63; 33 1039 81. JCT 3 121 67.**

$$\frac{-1 + 18z + (-(12z-1)^3)^{1/2}}{54 z^2}$$

1, 3, 12, 56, 288, 1584, 9152, 54912, 339456
**Réf. : CJM 15 269 63.**

$$\frac{3(1-8z)^{1/2} + 8z - 3(1-8z)^{3/2}}{4(1+(1-8z)^{1/2})^3 z}$$

1, 0, 4, 6, 24, 66, 214, 676, 2209, 7296, 24460, 82926, 284068, 981882, 3421318, 12007554, 42416488, 150718770, 538421590, 1932856590, 6969847484
**Réf. : CJM 15 265 63.**

$$\frac{(1+z)((-4z+1)^{3/2} - 1 + 6z - 6z^2 - 4z^3 - 6z^4) + 4z^5}{2(2z^5(z+2)^3(1+z))}$$

1, 3, 10, 33, 111, 379, 1312, 4596, 16266, 58082, 209010, 757259, 2760123, 10114131, 37239072, 137698584, 511140558, 1904038986, 7115422212, 26668376994
**Réf. : IC 16 351 70.**

$$\frac{1 - 3z - z^2 - (-(-1+4z)(-1+z+z^2)^2)^{1/2}}{2(2z^4 + z^5)}$$

1, 4, 15, 54, 193, 690, 2476, 8928, 32358, 117866, 431381, 1585842, 5853849, 21690378, 80650536, 300845232, 1125555054, 4222603968, 15881652606
**Réf. : IC 16 351 70.**

$$\frac{1 - 4z + z^2 + 2z^3 - (-(-1+4z)(z^2+2z-1)^2)^{1/2}}{2(2z^5 + z^6)}$$

1, 14, 120, 825, 5005, 28028, 148512, 755820, 3730650, 17978180, 84987760, 395482815
**Réf. : CAY 13 95. AEQ 18 385 78.**

$$\frac{1/2\,(1 - 21z + 180z^2 - 800z^3 + 1920z^4 - 2304z^5 + 1024z^6}{(z^5(4z-1)^5)}$$

$$-\;\frac{-(-(10z^4 - 50z^3 + 40z^2 - 11z + 1)^2(4z-1)^5)^{1/2})}{(z^5(4z-1)^5)}$$

1, 1, 1, 3, 16, 75, 309, 1183, 4360, 15783, 56750, 203929, 734722, 2658071, 9662093, 35292151, 129513736, 477376575, 1766738922, 6563071865, 24464169890
**Réf. : JCT B21 75 76.**

$$-\,\frac{1/2\,(-1 + 10z - 42z^2 + 98z^3 - 137z^4 + 112z^5 - 48z^6 + 8z^7}{(z^2(2z-1)^2(z-1)^4)}$$

$$+\;\frac{(-(-1+4z)(2z-1)^4(z-1)^8)^{1/2})}{(z^2(2z-1)^2(z-1)^4)}$$

1, 3, 9, 25, 69, 189, 518, 1422, 3915, 10813, 29964, 83304, 232323, 649845, 1822824, 5126520, 14453451, 40843521, 115668105, 328233969, 933206967, 2657946907, 7583013474
**Réf. : JCT A23 293 77.**

$$\frac{1 - 3z + 2z^3 - (-(3z^2 + 2z - 1)(-1+2z)^2)^{1/2}}{2z^6}$$

1, 4, 14, 44, 133, 392, 1140, 3288, 9438, 27016, 77220, 220584, 630084, 1800384, 5147328, 14727168, 42171849, 120870324, 346757334, 995742748, 2862099185
**Réf. : JCT A23 293 77.**

$$\frac{1 - 4z + 2z^2 + 4z^3 - z^4 - (-(-1+2z+3z^2)(1-3z+z^2+z^3)^2)^{1/2}}{z^8}$$

1, 5, 20, 70, 230, 726, 2235, 6765, 20240, 60060, 177177, 520455, 1524120, 4453320, 12991230, 37854954, 110218905, 320751445, 933149470, 2714401580, 7895719634
**Réf. : JCT A23 293 77.**

$$\frac{-\frac{1}{2}(-1 + 5z - 5z^2 - 5z^3 + 5z^4 + z^5)}{z^{10}} + \frac{(-(z+1)(3z-1)(z^2+z-1)^2(z^2-3z+1)^2)^{1/2}}{z^{10}}$$

1, 6, 27, 104, 369, 1242, 4037, 12804, 39897, 122694, 373581, 1128816, 3390582, 10136556, 30192102, 89662216, 265640691, 785509362, 2319218869, 6839057544
**Réf. : JCT A23 293 77.**

$$\frac{\frac{1}{2}(1 - 6z + 9z^2 + 4z^3 - 12z^4 + 2z^6)}{z^{12}} - \frac{(-(z+1)(3z-1)(z-1)^2(2z-1)^2(2z^2+2z-1)^2)^{1/2}}{z^{12}}$$

1, 2, 6, 16, 45, 126, 357, 1016, 2907, 8350, 24068, 69576, 201643, 585690, 1704510, 4969152, 14508939, 42422022, 124191258, 363985680, 1067892399, 3136046298, 9217554129
**Réf. : Comtet Louis, Advanced Combinatorics, p. 78.**

$$\frac{z + (z+1)^{1/2}(1-3z)^{1/2} - 1}{2(z^2(z+1)^{1/2}(1-3z)^{1/2})}$$

1, 3, 9, 26, 75, 216, 623, 1800, 5211, 15115, 43923
**Réf. : AAM 9 340 88.**

$$\frac{1 - 3z - (-(3z^2 + 2z - 1)(-1 + 2z)^2)^{1/2}}{2(3z^4 - z^3)}$$

# BIBLIOGRAPHIE

## Abbréviations des références